\documentclass{amsart}

\usepackage{amsfonts,amsthm,amsmath,amsfonts,latexsym,amssymb}

\usepackage{graphicx,upref}
\usepackage{color}
\usepackage{mathrsfs}
\usepackage[T1]{fontenc}
\usepackage[colorlinks=true,linkcolor=azul,citecolor=azul]{hyperref}

\newtheorem{fed}{Definition}[section]
\newtheorem{teo}[fed]{Theorem}
\newtheorem*{teo*}{Theorem}
\newtheorem{lem}[fed]{Lemma}
\newtheorem{cor}[fed]{Corollary}
\newtheorem{pro}[fed]{Proposition}

\theoremstyle{definition}

\newtheorem{rem}[fed]{Remark}

\newtheorem{exas}[fed]{Examples}

\def\bdem{\begin{proof}}
\def\edem{\end{proof}}
\def\ds{\displaystyle}

\def\bm{\left(\begin{array}}
\def\em{\end{array}\right)}
\def\ben{\begin{enumerate}}
\def\een{\end{enumerate}}
\def\barr{\begin{array}}
\def\earr{\end{array}}

\definecolor{azul}{rgb}{0.1,0.6,0.86}
\definecolor{titleblue}{rgb}{0.13,0.49,0.69}
\definecolor{mylred}{rgb}{0.85,0.24,0.2}
\definecolor{myblue}{rgb}{0,0.33,0.55}
\definecolor{myyellow}{rgb}{0.42,0.24,0.52}
\definecolor{mygreen}{rgb}{0.12,0.5,0.29}
\definecolor{myred}{rgb}{0.74,0.13,0.13}
\definecolor{mylblue}{rgb}{0.2,0.75,1}
\definecolor{mylgreen}{rgb}{0.68,0.98,0.6}
\definecolor{mylyellow}{rgb}{0.86,0.85,0.55}
\definecolor{myllyellow}{rgb}{0.87,0.86,0.56}
\definecolor{naranja}{RGB}{249,153,96}
\definecolor{sidebardarkcolor}{rgb}{0.21,0.31,0.40}
\definecolor{sidebarlightcolor}{rgb}{0.7,0.77,0.836}

\def\al{\alpha} 
\def\la{\lambda}
\def\eps{\varepsilon}
\def\fii{\varphi}
\def\w{\omega}
\def\z{\zeta}

\def\La{\Lambda}

\def\W{\Omega}

\def\N{\mathbb{N}}
\def\Z{\mathbb{Z}}
\def\R{\mathbb{R}}
\def\C{\mathbb{C}}

\def\D{\mathbb{D}}
\def\T{\mathbb{T}}

\def\1{\mathbb{1}}

\def\cA{\mathcal{A}}

\def\ele{\mathcal{L}}

\newcommand{\peso}[1]{ \quad \text{ \rm  #1 } \quad }

\newcommand{\sub}[2]{{#1}_{\mbox{\tiny{${#2}$}}}}

\DeclareMathOperator{\Preal}{Re}

\newcommand{\modu}[1]{\left|#1 \right|}
\newcommand {\lel}{\left\{}
\newcommand {\ril}{\right\}}
\newcommand{\norm}[1]{\left\|#1\right\|} 

\newcommand{\pint}[1]{\displaystyle \left \langle\, #1 \, \right\rangle}

\newcommand{\ol}{\overline}

\newcommand{\hil}{\mathcal{H}}

\newcommand{\conv}{\xrightarrow[n\rightarrow\infty]{}}

\newcommand{\hsd}{H^2(\D)}

\newcommand{\hst}{H^2(\T)}

\newcommand{\model}[1]{\hil(#1)}

\definecolor{mcolor}{rgb}{0,0.5,0.5}

\providecommand{\keywords}[1]{\textbf{\textbf{Keywords:}} #1}
\providecommand{\class}[1]{\textbf{\textbf{ AMS 2010 Mathematics Subject Classification:}} #1}

\begin{document}


\title[]{Model subspaces techniques to study Fourier expansions in $L^2$ spaces associated to singular measures}

\date{}

\author{Jorge Antezana}
\address{Departamento de Matemática, Universidad Nacional de La Plata and, 
 Instituto Argentino de Matemática "Alberto P. Calderón" (IAM-CONICET), Buenos Aires, Argentina}
\email{antezana@mate.unlp.edu.ar}
\thanks{The first author is partially supported by Grants CONICET-PIP 152, MTM2016-75196-P, PICT-2015-1505, and UNLP-11X681.}

\author{María Guadalupe García}
\address{Departamento de Matemática, Universidad Nacional de La Plata and, 
 Instituto Argentino de Matemática "Alberto P. Calderón" (IAM-CONICET), Buenos Aires, Argentina}
\email{mggarcia@mate.unlp.edu.ar}
\thanks{The second author is partially supported  by Grants CONICET-PIP 152, PICT-2015-1505 and UNLP-11X681.}

\maketitle

\begin{abstract}

\noindent Let $\mu$ be a probability measure on $\T$ that is singular with respect to the Haar measure. In this paper we study Fourier expansions in $L^2(\mathbb{T},\mu)$ using techniques from the theory of model subspaces of the Hardy space. Since the sequence of monomials $\{z^n\}_{n\in \N}$ is effective in $L^2(\mathbb{T},\mu)$, it has a Parseval frame associated via the Kaczmarz algorithm. Our first main goal is to identify the aforementioned frame with boundary values of the frame $P_\varphi(z^n)$ for the model subspace $\mathcal{H}(\varphi)= H^2 \ominus \varphi H^2$, where $P_\fii$ is the orthogonal projection from the Hardy space $H^2$ onto $\mathcal{H}(\varphi)$. The study of Fourier expansions in $L^2(\mathbb{T},\mu)$ also leads to consider positive kernels in the Hardy space. Our second main goal is to study the set of measures $\mu$ which reproduce a kernel contained in a model subspace. We completely characterize this set when the kernel is the reproducing kernel of a model subspace, and we study the consequences of this characterization.

 \medskip

 \noindent \keywords{Model subspaces; Kaczmarz algorithm; Fourier expansions; Parseval Frames}  \\
\noindent \class{ Primary: 42A16, 30H10 }
\end{abstract}
 
\tableofcontents

\section{Introduction}

Let $\{f_n\}_{n\geq 0}$ be a sequence of vectors  belonging to a (separable) Hilbert space $\hil$. This sequence is called a \textbf{frame} for $\hil$ if there exist constants $A,B>0$ such that for every $f\in\hil$ it holds that
$$
A\|f\|_{\hil}\leq \sum_{n\geq 0} \big|\pint{f,f_n}_{\hil}\big|^2\leq B\|f\|_{\hil}.
$$ 
If $A=B=1$, the frame is called \textbf{Parseval frame}. Given $\la\in\R^d$, let $e_\la$ denote the exponential function $e_\la(x)=e^{2\pi i (x\cdot \la)}$, where $x\cdot \la$ denotes the standard inner product of $\R^d$. 
If $\mu$ is a finite Borel measure on $\R^d$, we say that a family of exponentials $E(\La)=\{e_\la\}_{\la\in \La}$ is a \textbf{Fourier frame} for $L^2(\mu)$ if the family $E(\La)$ constitutes a frame for $L^2(\mu)$. This notion was introduced by Duffin and Schaeffer in \cite{DS}. 
It is known  that if a measure $L^2(\mu)$ admits a Fourier frame, then it must be of \textbf{pure type}, that is, the measure $\mu$ is either discrete, absolutely continuous or singular continuous (see for instance \cite{HLL, LW}).

\medskip

If $\mu$ is discrete, then it has a finite number of atoms. Therefore the analysis essentially is reduced to $\C^d$, where the Fourier frames are generator systems consisting on exponential functions. So, good characterization of such frames are obtained by using linear algebra techniques.

\medskip

If $\mu$ is absolutely continuous, then $\mu$ has to be supported on a set $\W$ of finite Lebesgue measure in $\R^d$, and its density function is bounded from above and below almost everywhere in $\W$ (see  \cite{DL}, \cite{Lai}). Moreover, if this is the case, then Nitzan, Olevskii and Ulanovskii proved in \cite{NOU} that $L^2(\mu)$ admits a Fourier frame. When $\W$ is bounded there are many ways to prove the existence of Fourier frames, and the proofs are not difficult with the techniques known nowadays. However, when $\W$ is unbounded, the existence of Fourier frames is much more difficult, and the proof in  \cite{NOU} uses some deep results in Functional Analysis such as the solution of the Kadison-Singer problem \cite{MSS}.

\medskip

The last case, when the measure $\mu$ is continuous singular, is much less understood. It is known that for many self-similar measures $\mu$ the corresponding space $L^2(\mu)$ admits a Fourier frame. Even in some cases they admit orthonormal basis of exponentials (see \cite{Dai}, \cite{DHL}, \cite{DLW}, \cite{DHS}, \cite{JP}). However, for instance, the question of whether or not the middle-third Cantor measure has or not a Fourier frame is still open. 

\medskip

Another possibility of Fourier expansions is in terms of the so called \textbf{effective sequences}.
Let  $\{\fii_n\}_{n \geq 0}$ be a complete sequence of unit vectors in a Hilbert space $\hil$.  Given $\xi\in\hil$, we define inductively
\begin{align}
\xi_0 &= \pint{\xi, \fii_0} \fii_0,    \nonumber\\
\xi_n &= \xi_{n-1} + \pint{\xi- \xi_{n-1}, \fii_n} \fii_n.  \label{Algoritmo-K}
\end{align}
If $\ds\lim_{n\to\infty}\|\xi-\xi_n\| = 0$ for every $\xi\in\hil$, then the sequence $\{\fii_n\}_{n \geq 0}$ is called effective. The above inductive construction  that leads to the vectors $\xi_n$ is the so called Kaczmarz algorithm. This algorithm was introduced for finite dimensional spaces by  Kaczmarz \cite{K}, and it was studied in infinite dimensional Hilbert space by Kwapie$\acute{\text{n}}$ and Mycielski in \cite{KM} (see also \cite{H-W}). The vectors $\xi_n$ obtained by the Kaczmarz algorithm can be expressed as
\begin{equation}
\xi_{n} = \ds{\sum_{k=0}^{n} \pint{\xi, \gamma_k} \fii_k}, \label{eq intro xn con los gn} 
\end{equation}
where the vectors $\gamma_n$ are also defined by a recurrence relation: 
\begin{align}
\gamma_0 &= \fii_0\nonumber \\
\gamma_n &= \fii_n - \ds{\sum_{k=0}^{n-1}} \pint{\fii_n, \fii_k} \gamma_k. \label{eq intro recurrence for the gn}
\end{align}
The important fact is that, $\{\fii_n\}_{n\geq 0}$ is effective if and only if $\lel \gamma_n \ril_{n\geq 0}$ is a Parseval frame. In general, an effective sequence $\{\fii_n\}_{n\geq 0}$ is (super) redundant, and there are different ways to reconstruct a vector $\xi$. In particular, there may be more than one Parseval frame that satisfy \eqref{eq intro xn con los gn}. The advantage of $\lel \gamma_n \ril_{n\geq 0}$ is its recursive definition in terms of the $\fii_n$. 

\medskip

In \cite{KM} (see also \cite{H-S}), the following result was proved.

\begin{teo}\label{efectiva + medida}
Let $\mu$ be a Borel probability measure in the interval $[0,1)$. The sequence of integer exponentials 
$\lel e_n\ril_{n\geq 0}$ is effective if and only if  $\mu$ is either the normalized Lebesgue measure restricted to $[0,1)$, or singular with respect to Lebesgue measure. 
\end{teo}

Note that, given a singular measure $\mu$, which can be any singular Cantor measure, for every $f\in L^2(\mu)$ we have the Fourier expansion 
$$
f=\sum_{n=0}^\infty \pint{f,g_n}e^{2\pi i n t},
$$ 
where $\{g_n\}_{n\geq 0}$ is a Parseval frame. Since the effective sequence is super-redundant, it is not a Bessel sequence. So, the relevance of these Fourier expansions lie in the Kaczmarz algorithm that is behind them, and in the Parseval frame $\{g_n\}_{n\geq 0}$.

\medskip

The main aim of this paper is to use some tools from the theory of model subspaces of the Hardy space $H^2$ in the study of Fourier representation for $L^2$ spaces of singular measures (see Section 2 for basic definitions on Hardy spaces and model subspaces). This idea of course is not new (see for instance \cite{HNP}). In the case of effective sequences, the Hardy space techniques appear in the work of Haller and Szwarc \cite{H-S}. The theory of model spaces is also present in the work by Herr, Jorgensen and Weber \cite{H-J-W}.

\medskip

Let $\mu$ be a Borel measure in $\T$ that is singular with respect to the Lebesgue measure. Our first main result is a representation of the Parseval frame associated to the effective sequence of monomials in $L^2(\mu)$. We characterize this Parseval frame as the boundary values of the Parseval frame obtained by projecting the monomials onto a convenient model space (see Theorem \ref{proyecciones + frame}). We also analyze some consequences of this representation.   

\medskip

On the other hand, Dutkay and Jorgensen studied in \cite{DorinPalle} the connection between lacunary Fourier expansions in $L^2$ spaces associated to fractal measures and closed subspaces of the Hardy space. This led them to the study of positive matrices in the Hardy space with prescribed boundary representations, a study that was continued by Herr et. al. in \cite{H-J-W}. In this direction,  using some results on isometric embedding of model spaces due to Aleksandrov, we study the set of measures that reproduce a kernel contained in a model subspace 
(see Subsection 3.1 for the precise definitions).  When the kernel coincides with the reproducing kernel of a model subspace, we prove that those measures coincide with the set of measures $\mu$ such that the corresponding model subspace can be isometrically embedded in 
$L^2(\mathbb{T},\mu)$ (see Subsection \ref{subsection 32} for more details).

\bigskip

The paper is organized as follows. In section 2 we recall some definitions and results on Hardy spaces and model subspaces that we will need in the sequel. Section 3 contains the main results. Firstly, we will prove the aforementioned representation of the Parseval frame associated to the effective sequence of monomials. The second part of this section is devoted to the  study of positive kernels in the Hardy space.

\section{Preliminaries on Model subspaces}

In this section we recall some definitions and results related with model subspaces that we will need throughout this paper. The section is based on the monographs  \cite{C-M-R}, \cite{Nik1}, \cite{Nik2}, \cite{Rudin} and \cite{Sarason}.

\bigskip

The \textbf{Hardy space $\hsd$} is defined as the set of holomorphic functions $f:\D\to \C$ such that
$$
f(z)=\sum_{n=0}^\infty\, a_n z^n \quad  \mbox{and} \quad  \{a_n\}\in\ell^2(\N_{0}).
$$
This space is endowed with the norm 
\begin{equation}\label{eq norma hsd 2}
\norm{f}_2^2 := \sum_{n=0}^{\infty}\, \modu{a_{n}}^{2}.
\end{equation}
This norm provides a Hilbert space structure on $\hsd$. With respect to this Hilbert space structure, the evaluations are bounded functionals. Hence $H^2(\D)$ turns out to be a reproducing kernel Hilbert space, whose reproducing kernel is the so called \textbf{Szegö kernel} given by 
$$
k_{w}(z)= \frac{1}{1-z\overline{w}},\quad w \in \D.
$$

Let $\T=\{z\in\C: \ |z|=1\}$, and let $m$ denote the Lebesgue measure normalized such that its total mass is one. Using this measure, the norm in $\hsd$ has the following alternative expresion
\begin{equation}\label{eq norma hsd 3}
\|f\|_2^2=\lim_{r\to1^-}\frac{1}{2\pi}\int_\T |f(r\z)|^2 dm(\z).
\end{equation}
Given $f\in \hsd$, for $m$-almost every $\w\in\T$ there exists the non-tangential limit (see \cite{Rudin} for its definition)
$$
F(\w) ={\lim_{\angle z\rightarrow \w}} f(z).
$$
It is not difficult to prove that $F$ defines an element of $L^{2}(\T,m)$ that satisfies the condition
\begin{equation}\label{eq HST}
\int_\T F(\z) \overline{\z^n} \ d\z= 0,\quad   \forall\, n<0.
\end{equation}
The subspace consisting of those functions of $L^{2}(\T,m)$ that satisfy \eqref{eq HST} is called \textbf{Hardy space of the torus}, and we will denote it by $\hst$. The map $f\mapsto F$ using non-tangential limits is an isometric isomorphism between the spaces $\hsd$ and $\hst$. The inverse of this map is given by the Poisson integral. Recall that given a finite Borel measure $\mu$ in $\T$, the \textbf{Poisson integral} of $\mu$ is defined by
$$
P[\mu](z)=\int_\T  P_z(\z) \, d\mu(\z), 
$$ 
where $\ds P_z(\z)=\Preal\left(\frac{\z + z}{\z - z}\right)=\frac{1-\modu{z}^2}{\modu{\z - z}^2}$ is the well known Poisson kernel of the disc. So, given $F\in\hst$, the function
$$
f(z)=P[Fdm](z)=\int_\T\frac{1-\modu{z}^2}{\modu{\z - z}^2}\,F(\z)\, dm(\z) 
$$
belongs to $\hsd$, and $\ds F(\z) ={\lim_{\angle z\rightarrow \z}} f(z)$. Moreover, $\|f\|_2=\|F\|_2$, where the second norm is the $L^2$ norm  with respect to the measure $m$. 

\paragraph{Convention:} From now on, many times we will identify the spaces $\hsd$ and $\hst$, using one or the other depending on our convenience. In this case we will simply write $H^2$. Moreover, we will use the same letter for an element $F\in \hst$ and its Poisson extension to the disc, which as we mentioned is the element in $\hsd$ corresponding to $F$ by the aforementioned isometric isomorphism.

\bigskip
Let $S$ denote the well known shift operator, defined for $f\in \hsd$ by
$$
Sf(z)=zf(z).
$$
Its definition in $\hst$ is given through the aforementioned isomorphism. One of the most important results in the Hardy spaces theory is the Beurling's theorem, which gives a precise description of Shift invariant subspaces. 

\begin{teo}[Beurling]\label{teo Beurling}
Let $\ele$ be a non-zero shift invariant subspace of $\hst$, then there is a function  
$\vartheta\in \hst$ such that 
\begin{equation}\label{eq Beurling theorem}
|\vartheta(\w)|=1,
\end{equation}
for almost every $\w\in\T$, and 
$$
\ele=\vartheta\hsd:=\{\vartheta f:\ f\in\hsd\}.
$$ 
Reciprocally, if $\vartheta\in \hst$ and satisfies \eqref{eq Beurling theorem} for almost every $\w\in\T$, then $\vartheta H^2$ is a shift invariant subspace of $\hst$.
\end{teo}

A function $\vartheta\in \hst$  that satisfies \eqref{eq Beurling theorem}, for almost every $\w\in\T$, is called \textbf{inner function} or simply \textbf{inner}. It is not difficult to prove that the modulus of an inner function, as an element of $\hsd$, is bounded by one. Moreover, given two inner functions $\vartheta_1$ and $\vartheta_2$, we say that \textbf{$\vartheta_1$ divides $\vartheta_2$} if $\vartheta_2/\vartheta_1\in\hsd$. Equivalently  $\vartheta_1$ divides $\vartheta_2$ if $\vartheta_1H^2$ contains the subspace $\vartheta_2H^2$.

\medskip

\begin{fed}
Let $\fii$ be an inner function. The space $\model{\fii}:= H^2 \ominus \fii H^2$ is called \textbf{model space}.    
\end{fed}

\medskip

The orthogonal projection $P_\fii$ from $L^2(\T, m)$ onto the model space $\model{\fii}$ is defined by 
\begin{equation}\label{eq proyeccion modelo 1}
P_\fii f(\w) = f(\w) - \fii(\w) \sub{P}{H^2}(\overline{\fii}(\w) f(\w)), 
\end{equation}
where $\sub{P}{H^2}: L^2(\T, m) \rightarrow H^2(\T)$ is the projection onto the Hardy space given by
$$
\sub{P}{H^2} f(\w) = \lim_{\angle z \rightarrow \w} \int_\T \frac{f(\w)}{1- z \overline{\w}} \, d\w.
$$

Clearly, the model subspaces (as subspaces of $\hsd$) are also reproducing kernel Hilbert spaces. 
The reproducing kernels $k_z^\fii$ for the model space $\model{\fii}$ take the form
$$
k_z^\fii (w) = \frac{1 - \overline{\fii(z)} \fii(w)}{1 -\overline{z} w},
$$
with $z,w \in \D$. Using this reproducing kernel we have the following alternative form to write the orthogonal projection onto $\model{\fii}$ 

\begin{equation}\label{eq proyeccion modelo 2}
(P_\fii f) (z) 
= \pint{f, k_z^\fii} 
= \int_\T f(\z) \, \frac{1-\fii(z) \overline{\fii(\z)}}{1-z\overline{\z}} \,dm(\z).
\end{equation}

\bigskip

The following result provides a dense subspace of $\model{\fii}$, which will be very useful in the sequel.

\begin{teo} \label{densidad en el esp. modelo}  
Let $A$ be the disk algebra, i.e $A$ is the space of continuous functions on $\overline{\D}$ which
are analytic on $\D$. If $\fii$ is an inner function then
$$
\cA(\fii) := A \cap \model{\fii}
$$
is a dense subset of $\model{\fii}$.
\end{teo}

\subsection{Clark measures}

Let $\fii:\D \rightarrow \D$ be an analytic function.  Given $\al \in \T$, the function  
$$
\Preal\left(\frac{\al+\fii(z)}{\al -\fii(z)}\right) = \frac{1-\modu{\fii(z)}^{2}}{\modu{\al-\fii(z)}^{2}}
$$  
is harmonic and positive. Then, by Herglotz's Theorem, there is a unique finite positive measure 
$\mu_{\al}$ on $\T$ such that 
$$
\Preal\left(\frac{\al + \fii(z)}{\al-\fii(z)}\right) = \int_{\T} P_{z}(\z)\,d\mu_{\al}(\z).
$$ 
The measures $\mu_{\al}$, with $\al \in \T$, are called \textbf{Clark measures} corresponding to $\fii$.
If $\fii$ is an inner function, each measure $\mu_{\al}$ is singular with respect to the Lebesgue measure. Moreover
$$
\mu_\alpha \big(\T\setminus \big\{\w\in\T:\ \lim_{\angle z \rightarrow \w}\fii(z)=\alpha\big\}\big)=0.
$$

\medskip

These measures are named after a work of D.N. Clark, who in \cite{Clark} introduced them in his study of spectral measures of unitary extension of the restriction of the shift operator to model subspaces (see also \cite{C-M-R}, \cite{Eero}, \cite{Sarason} and the references therein).

\medskip

Clark showed that the operator $V_\al$ defined on $L^2(\T, \mu_\al)$ by 
$$
\left(V_\al g\right)(z) := \left( 1- \overline{\al} \fii(z)\right) \int_\T \frac{g(\z)}{1 - \overline{\z}z} \,d\mu_{\al}(\z)
$$ 
is an isometric isomorphism from $L^2(\T, \mu_\al)$ onto the model space $\model{\fii}$. 
The inverse of $V_\alpha$ was characterized by Poltoratskii in \cite{Poltoratskii}. More precisely he proved the following theorem, which will be very important for this work.

\begin{teo}[Poltoratskii]\label{Polto}
Let $\alpha\in\T$, and $f\in \model{\fii}$. Then, for $\mu_\alpha$-almost every $\w\in\T$ the non-tangential limits
$$
\lim_{\angle z \rightarrow \w}f(z)
$$
exist. Moreover, if $f^*$ denotes the element of $L^2(\T, \mu_\alpha)$ defined by these limits, then
$$
\|f\|_{\model{\fii}}=\|f\|_{H^2}=\|f^*\|_{L^2(\T, \mu_\alpha)},
$$
and $V_\alpha(f^*)=f$.
\end{teo}

\begin{rem}
Many times throughout the paper we will use the notation $\mu^{\fii}_\alpha$ to emphasize this Clark measure is related to the function $\fii:\D\to\D$.
\end{rem}

\section{Main results}

\subsection{Parseval frames associated to effective sequences}

Let $\mu$ be a Borel probability measure on $\T$ that is singular with respect to the Lebesgue measure on the torus, and identify  the interval $[0,1)$ with $\T$ in the usual way. According to the result of  Kwapie$\acute{\text{n}}$ and Mycielski, the monomials $\{\w^n\}_{n\geq 0}$ constitute an effective sequence for $L^2(\T, \mu)$. As a consequence, every element $f\in L^2(\T,\mu)$ can be written as
\begin{equation}\label{eq main result effectivas}
f=\sum_{n=0}^\infty \pint{f,g_n} \w^n
\end{equation}
where $\lel g_n\ril_{n\geq 0}$ is a Parseval frame for $L^2(\T,\mu)$. The convergence of the series is understood in the $L^2$-sense.

\medskip

Let $\psi:\D\to\C$ be the holomorphic function defined by the Herglotz transform as
$$
\psi(z)=\int_\T \frac{\z + z}{\z - z} \, d\mu(\z).
$$ 
Note that the real part of $\psi(z)$ is precisely $P[\mu](z)$, which is positive. Therefore, the function
\begin{equation}\label{eq la fi}
\fii(z)=\frac{\psi(z)-1}{\psi(z)+1} 
\end{equation}
is a self-map of the disc. A little algebra shows that
$$
P[\mu](z)=\frac{1-|\fii(z)|^2}{|1-\fii(z)|^2}.
$$
Since $\mu$ is singular with respect to the Lebesgue measure on $\T$, by a standard result on Poisson integral (see for instance 
\cite[Thm 11.22]{Rudin}) it holds that for $m$-almost all $\w\in\T$
$$
\lim_{\angle z\rightarrow \w} P[\mu](z)=0,\quad\mbox{which implies that}\quad \lim_{\angle z\rightarrow \w} |\fii(z)|=1,
$$
for $m$-almost all $\w\in\T$. In consequence $\fii$ is an inner function. Moreover, since $\mu$ is a probability measure, $\psi(0)=1$ and therefore $\fii(0)=0$. 

\medskip

Let $\model{\fii}$ be the associated model space. As the monomials $\{z^n\}_{n\geq 0}$ constitute an orthonormal basis for the Hardy space, their projection $\{P_\fii(z^n)\}_{n\geq 0}$ to the model space $\model{\fii}$ is a Parseval frame. 

\bigskip

The following theorem is one of the main results of this paper. It relates this Parseval frame $\{P_\fii(z^n)\}_{n\geq 0}$ with the Parseval frame $\lel g_n\ril_{n\geq 0}$ associated to the monomial when they are considered as an effective sequence for $L^2(\T,\mu)$.

\medskip

\begin{teo}\label{proyecciones + frame}
For $\mu$-almost every $\z\in\T$
$$
g_n(\z)=\lim_{\angle w \rightarrow \z}P_\fii(z^n)(w).
$$
\end{teo}

\bdem
Given $n\in\Z$, consider the Fourier coefficients of $\mu$  
$$
\widehat{\mu}(n)= \int_\T \overline{\w}^{\,n} d\mu(\w).
$$
Define $M$ as the lower triangular Toeplitz matrix whose entries are 
$$
M_{ij}=\begin{cases}
\widehat{\mu}(i-j)&\mbox{if $i> j$}\\
0&\mbox{if $i\leq j$}
\end{cases}.
$$
Let $U$ the lower triangular Toeplitz matrix defined by the algebraic relation
\begin{equation}\label{eq matrices de Toeplitz}
(M+I)^{-1}=I+U,
\end{equation}
where $I$ denotes the identity matrix. Let $\{u_n\}_{n\geq 1}$ such that 
$$
U_{ij}=\begin{cases}
u_{i-j}&\mbox{if $i> j$}\\
0&\mbox{if $i\leq j$}
\end{cases}.
$$
Let $\phi$ and $\eta$ be the functions defined by
\begin{equation}
\phi(z) = \sum_{n= 1}^\infty -\overline{u_n} z^n\peso{and} \eta(z)=\sum_{n=1}^\infty \widehat{\mu}(n) z^n. \label{eq: serie de inner}
\end{equation}
Note that the identity \eqref{eq matrices de Toeplitz} implies that $(1+\eta(z))(1+\phi(z))=1$. From this identity we get that
$$
1+\eta(z)+\overline{\eta(z)}=\frac{1-|\phi(z)|^2}{|1-\phi(z)|^2}.
$$
This fact was observed by Haller and Szwarc, who used it to prove that $\phi$ is inner (see Section 4 in \cite{H-S}). However, if we also observe that
$$
1+\eta(z)=\int_{\T} \frac{1}{1-\overline{\z}z}\, d\mu(\z)
$$
we easily obtain that 
$$
\frac{1-|\phi(z)|^2}{|1-\phi(z)|^2}=P[\mu](z).
$$
This shows that $\phi$ and the function $\fii$ defined by (\ref{eq la fi}) have the same real part. Since $\phi(0)=\fii(0)=0$, we conclude that  $\phi(z)=\fii(z)$. Using this, if 
$w \in \D$

\begin{align*}
\left(P_\fii z^n\right) (w) 
&= w^n -\fii(w) \int_\T \frac{\overline{\fii(\z)} \, \z^n}{1- w \overline{\z}}\, d\z
= w^n -\fii(w) \int_\T \overline{\fii(\z)}\,\z^n \ds{\sum_{k=0}^\infty} \left(w \overline{\z}\right)^k d\z \\
&= w^n -\fii(w) \ds{\sum_{k=0}^\infty} \underbrace{\left(\int_\T \overline{\fii(\z)}\,\z^{n-k} 
d\z \right)}_{=0 \, \text{if}\, k\geq n} w^k 
\\&= w^n -\fii(w) \ds{\sum_{k=0}^{n-1}}\left(\int_\T \overline{\fii(\z)}\,\z^{n-k} 
d\z \right) w^k \\
&= w^n -\fii(w) \ds{\sum_{k=0}^{n-1}} -\, u_{n-k} \, w^k 
= w^n + \ds{\sum_{k=0}^{n-1}} \fii(w) \, u_{n-k}\, w^k .\\
\end{align*}
Now, if $w$ tends to $\z$ non-tangentially, we get
\begin{align*}
(P_\fii z^n)^* (\z) 
&= \z^n + \sum_{k=0}^{n-1} \, \fii(\z) \, u_{n-k}\, \z^k 
= \z^n + \sum_{k=0}^{n-1}\, u_{n-k}\, \z^k, \\ 
\end{align*}
where the second equality is because $\fii =1 $ on the support of $\mu$. It only remains to show that the last expression is precisely 
$g_n(\z)$. With this aim, note that as a consequence of \eqref{eq intro recurrence for the gn}
$$
\z^n = \sum_{j=0}^n \widehat{\mu}(n-j) g_j(\z), 
$$
which in terms of the matrix $M$ defined above can be written as
$$(I+M) 
\begin{pmatrix}
         g_0     \\ 
         g_1      \\
         \vdots  \\
          g_n    \\
					\vdots \\
\end{pmatrix}
=\begin{pmatrix}
         1         &  0         &    0     &  \cdots  \\ 
         m_{1,0}   &  1         &    0     &  \cdots  \\
         m_{2,0}   &  m_{2,1}   &    1     &  \cdots   \\
          \vdots   &  \vdots    &  \vdots  &  \vdots  \\
         m_{n,0}   &  m_{n,1}   &  m_{n,2} &  \cdots   \\ 
          \vdots   &   \vdots   &  \vdots  &   \vdots  \\
\end{pmatrix}
\begin{pmatrix}
         g_0     \\ 
         g_1      \\
         \vdots  \\
          g_n    \\
	\vdots \\
\end{pmatrix}
=\begin{pmatrix}
         \w_0     \\ 
         \w_1      \\
         \vdots  \\
         \w_n    \\
	\vdots \\
\end{pmatrix}.
$$
From this we conclude that $g_n(\z)= \z^n + \ds{\sum_{k=0}^{n-1}} \, u_{n-k} \, \z^k$, which concludes the proof.
\edem

Note that $\mu$ is the Clark measure associated to $\fii$ and $\alpha=1$. Since $\fii$ is inner, the other Clark measures $\mu_\alpha^{\fii}$ are also singular with respect to the Lebesgue measure on $\T$. So, the monomials are also effective in $L^2(\T,\mu_\alpha)$, and they have the correspondig Parseval frames $\{g_n^{(\alpha)}\}_{n\geq 0}$ associated. Since the Clark measure $\mu_\alpha^{\fii}$ is the Clark measure $\mu_1^{(\overline{\alpha}\fii)}$ related to the inner function $\overline{\alpha}\fii$, from the previous result we directly get the next corollary.

\begin{cor}\label{cor parrafo}
Given $\alpha\in\T$, for $\mu_\alpha^{\fii}$-almost every $\z\in\T$
$$
g_n^{(\alpha)}(\z)=\lim_{\angle w\rightarrow \z}P_\fii(z^n)(w).
$$
\end{cor}

\bigskip

Before going on, we would like to point out that Theorem \ref{proyecciones + frame} has the following simple consequence: for every $f\in\model{\fii}$ by Poltoratskii's theorem
\begin{align}
f(z) &= \ds{\sum_{n=0}^\infty} \pint{f, z_n}_{H^2} z^n
= \ds{\sum_{n=0}^\infty} \pint{P_\fii f, z_n}_{\model{\fii}} z^n \nonumber\\
&= \ds{\sum_{n=0}^\infty} \pint{f, P_\fii z_n}_{\model{\fii}} z^n 
= \ds{\sum_{n=0}^\infty} \pint{f^*, g_n}_{L^2(\mu)} z^n, \label{eq the key1}
\end{align}
where $f^*$ is the element of $L^2(\T,\mu)$ defined by taking non-tangential limits. Note that, since the sequence of monomials is overcomplete, this simple argument can not be made without the knowledge of Theorem \ref{proyecciones + frame}.

\medskip

Now, as an application of  Theorem \ref{proyecciones + frame}, we will provide simpler and more conceptual proof of a result appeared in 
\cite{H-W} related with the so called \textbf{normalized Cauchy transform}. Given a Borel (complex) measure on $\T$ we define the Cauchy transform of $\mu$ as
$$
K(\mu)(z)=\int_\T \frac{1}{1-z\overline{\z}}\,d\z,
$$
where $z\in\D$. Since $\fii(0) = 0 $, using the formula
$$
(K \,d\mu_\al^\fii)(z) = \frac{1}{1-\ol{\al}\fii(z)},
$$  
the operator $V_\al$ introduced in the preliminaries can be written as
$$
(V_\al g) (z) = \frac{K(g \,d\mu_\al^\fii)(z)}{(K \,d\mu_\al^\fii)(z)}.
$$
The right hand side of this identity is called \textbf{normalized Cauchy transform}.  

\begin{cor}\label{nucleo de Cauchy + frame}
Let $\mu$ a Borel probability measure that is singular with respect to the Lebesgue measure on $\T$. If $\lel g_n \ril_{n\in \N}$ is the Parseval frame for
$L^2(\T,\mu)$ associated to the effective sequence of the monomials, then 
$$
\frac{K(g \,d\mu)(z)}{(K \,d\mu)(z)}=\ds{\sum_{n=0}^\infty} \pint{g, g_n}_{L^2(\mu)} z^n.
$$
\end{cor}
\bdem
Let $\fii$ be the inner function defined by \eqref{eq la fi}. Since the normalized Cauchy transform coincides with the isometric isomorphism $V_1$ between $\model{\fii}$ and $L^2(\T,\mu)$, using \eqref{eq the key1} we get for every $g \in L^2(\T,\mu)$,
$$
V_1 g 
= \ds{\sum_{n=0}^\infty} \pint{V_1^* g, g_n}_{L^2(\mu)} z_n
= \ds{\sum_{n=0}^\infty} \pint{g, g_n}_{L^2(\mu)} z^n, 
$$ 
where in the second equality we have used Poltoratskii's theorem (Theorem \ref{Polto} of the preliminary section). 
\edem

\medskip 

Another application is the following expression for the reproducing kernel of $\model{\fii}$.

\begin{cor}
Let $\mu$ be a Borel probability measure that is singular with respect to the Lebesgue measure on $\T$, and $\fii$ the inner function defined by \eqref{eq la fi}. If $\lel g_n \ril_{n\in \N}$ is the Parseval frame for
$L^2(\T,\mu)$ associated to the effective sequence of the monomials, then 
$$
k_z^\fii (w) = \ds{\sum_{m=0}^\infty  \sum_{n=0}^\infty} \pint{g_n, g_m}_{L^2(\mu)}\overline{z}^n w^m.
$$
\end{cor}

\bdem
Indeed, using that the monomials form an orthonormal basis of the Hardy space we get
\begin{align*}
k_z^\fii (w) &=\sum_{m=0}^\infty \pint{k_z^\fii , w_n}_{H^2} w^n = \sum_{m=0}^\infty \pint{k_z^\fii ,P_\fii(w_n)}_{\model{\fii}}w^n\\
&=\sum_{m=0}^\infty \overline{P_\fii(w^n)(z)} w^n=\sum_{m=0}^\infty\sum_{n=0}^\infty \pint{g_n,g_m}_{L^2(\mu)} \overline{z}^n w^n,
\end{align*}
where in the last identity we have used \eqref{eq the key1} and Theorem \ref{proyecciones + frame}.
\edem

\subsection{Positive matrices in the Hardy space} \label{subsection 32}

Given a function $f:\D\to\C$, for any $r\in (0,1)$ by means of $f_r$ we will denote the function defined on $\T$ by
$$
f_r(\w)=f(r\w).
$$ 
Let $\mu$ be a finite positive Borel measure on $\T$. Throughout this section we say that $f^*\in L^2(\T,\mu)$ is an extension of $f$ if
\begin{equation}
\lim_{r\to 1^-} \|f^*-f_r\|_{L^2(\T,\mu)}=0. \label{eq: extension en L2}
\end{equation}
It is well know that any $f\in\hsd$ admits an extension with respect to the Lebesgue measure, and $f^*$ can be identified by taking non-tangential limits.

\medskip

Now, let us consider a positive definite kernel $k:\D\times\D\to\C$. As usual, by means of $k_z$ we will denote the function defined on $\D$ as $k_z(\cdot)=k(\cdot,z)$. Let $\mu$ be a positive Borel measure on $\T$, and suppose that for every $z\in\D$ the functions $k_z$ admit an extension to $L^2(\T,\mu)$. We say that $k$ reproduces itself with respect to $\mu$ if 
$$
k_z(w) = \int\limits_\T k_z^*(\z) \, \overline{k_ w^*(\z)} \, d\mu(\z).
$$
Given a fix kernel $k$, then by means of $M(k)$ we denote the set of all the positive Borel measure on $\T$ with respect to which $k$ reproduces itself. Finally, we will say that a positive kernel $k$ is included in $\model{\fii}$ if $k_z \in \model{\fii}$, for all $z\in \D$.

\begin{exas}
$\, $
\begin{itemize} 
\item The  Szegö kernel reproduces itself with respect to the Lebesgue measure.
\item Let $\fii$ an inner function, and let $k^\fii$ be the reproducing kernel of the model subspace $\model{\fii}$. Then, we will see that $k^\fii$ not only reproduces itself with respect to the Lebesgue measure, but it also reproduces itself with respect to any Clark measure $\mu_\alpha^\fii$.
\end{itemize}
\end{exas}

\bigskip

The first main goal of this section will be to relate the set $M(k^\fii)$, corresponding to the reproducing kernel of a model subspace $\model{\fii}$, with the following set $\text{Iso}(\fii)$ introduced by  Aleksandrov in \cite{Alek}.

\begin{fed}[Aleksandrov \cite{Alek}]
Given an inner function $\fii$,  the set $\text{Iso}(\fii)$ consists of those positive measures $\mu$ on $\T$ such that for every $g \in \cA(\fii)$ it holds that
\begin{equation}
\int\limits_\T \modu{g^*}^2 d\mu = \norm{g}^2_{H^2}, \label{eq: Iso}
\end{equation}
where $g^*$ denote the restriction of $g$ to $\T$.
\end{fed}

\medskip
If $g\in\cA(\fii)$, note that the functions $g_r$ converge uniformly to the restriction $g^*$. Since $\mu$ is finite, the uniform convergence implies that 
$$
\lim_{r\to 1^-} \|g^*-g_r\|_{L^2(\T, \mu)}=0.
$$
Therefore, if $g\in\cA(\fii)$ the notation $g^*$ for the restriction to the boundary is consistent with the notation introduced in \eqref{eq: extension en L2}. As this subspace is dense in $\model{\fii}$, equation \eqref{eq: Iso} defines an isometry $\Phi$ from $\model{\fii}$ into $L^2(\T,\mu)$. A priori the action on an element $f\in\model{\fii}\setminus \cA(\fii)$ is formal. However, the following result of Aleksandrov shows that this operator acts on each $f\in\model{\fii}$ by taking non-tangential limits.

\begin{teo}[Aleksandrov \cite{Alek2}]\label{Teo Alek}
If $\mu\in \text{Iso}(\fii)$, then for every $f\in \model{\fii}$
$$
\Phi(f)(\w)=\lim_{\angle z\rightarrow \w}f(z),
$$
holds for $\mu$-almost every $\w\in \T$.
\end{teo}

Our first main result is the following.

\begin{teo}\label{eq. con Iso}
Let $\fii$ be an inner function and let $k^\fii$ be the reproducing kernel for the model space 
$\model{\fii}$. Then $M(k^\fii)=\text{Iso}(\fii)$. Even more, for every $f\in \model{\fii}$
$$
\lim_{r\to 1^-} \|\Phi(f)-f_r\|_{L^2(\T,\mu)}=0.
$$
 
\end{teo}

\medskip

Since the linear span of the reproducing kernels $k^\fii_z$ is dense in $\model{\fii}$, this result is quite natural. The main difficulty in the proof of the theorem \ref{eq. con Iso} is to deal with the difference between the notions of boundary value used in the definitions of the sets $M(k^\fii)$ and $\text{Iso}(\fii)$. On one hand, in the definition of $M(k)$ the extension of the reproducing kernels is defined by means of $L^2(\T,\mu)$ limits of the radial restrictions (see \eqref{eq: extension en L2} above). On the other hand,  in the case of $ \text{Iso}(\fii)$, as we observed in the aforementioned Theorem of Aleksandrov, the definition of the boundary values is by means of non-tangential limits. 

\medskip

Since the proof of this result is quite long and technical, we leave it to the end of this section. Now, we continue with some consequences of Theorem \ref{eq. con Iso}. First of all, by the Poltoratskii's theorem and the definition of $\text{Iso}(\fii)$ we have the following result.   
\begin{cor}
All the Clark measures $\mu_\alpha^\fii$ belong to $M(k^\fii)$.
\end{cor}
On the other hand, in \cite{H-J-W}, Herr et. al. proved that if $k$ is a positive kernel included in a model subspace$ \model{\fii}$, and the Clark measure $\mu_1^\fii \in M(k)$ then 
$$
\mu_1^{\fii^n} \in M(k),
$$
for all $n\geq 0$.  The following proposition gives a generalization of this result.

\begin{pro}\label{medidas en M(K)}
Let $\fii$ be an inner function and let $k$ be a positive kernel included in $\model{\fii}$. If
$\mu_1^\fii \in M(k)$ and $\psi$ is an inner function such that $\fii | \psi$ then 
$\mu_\alpha^\psi \in M(k)$ for every $\alpha\in\T$.
\end{pro}

\bdem
It is enough to prove the result for $\alpha=1$ (see the paragraph before Corollary \ref{cor parrafo}). If $\psi$ is an inner function such that $\fii | \psi$ then
$\model{\fii} \subseteq \model{\psi} $.  Thus $k \subseteq \model{\psi} $.
Therefore there exists $ k^*_z \in L^2(\mu_1^\fii)\cap L^2(\mu_1^\psi)$ (defined by non-tangential limits  as well as by $L^2$ limit of the radial restrictions) such that 
\begin{align*}
k_z(w) &= \pint{k^*_z, k^*_w}_{L^2(\mu_1^\fii)} 
= \pint{k_z,k_w}_{\model{\fii}}
= \pint{k_z,k_w}_{H^2} \\
&= \pint{k_z,k_w}_{\model{\psi}}
= \pint{k^*_z,k^*_w}_{L^2(\mu_1^\psi)} 
= \int_\T k_z^*(\z) \overline{k_w^* (\z)} d \mu_1^\psi(\z),
\end{align*}
where we obtain the first equality because $\mu^\fii_1 \in M(k)$, the second and the fifth ones by the isometry $V_1$, the third and the fourth ones because the norms in $H^2$ and $\model{\fii}$ are the same. 
\edem

\medskip
Note that we can also write $\fii^n$ in the result of Herr et. al. as the composition $z^n \circ \fii$. It is a classical result that the composition of inner functions is an inner function (see for instance \cite[p. 323]{Steph}). So, as a consequence of the previous proposition we get the following variation on the result of Herr, Weber and Jorgensen aforementioned.

\medskip

\begin{cor}\label{Cor inner inner}
Let $\vartheta$ be an inner function such that $\vartheta(0)=0$ and let $k$ be a positive kernel in $H^2$ which is included in $\model{\fii}$. If $\mu_1^\fii \in M(k)$ then the function $\psi = \vartheta \circ \fii$ is inner and $\mu^\psi_\alpha \in M(k)$ for every $\alpha\in\T$.
\end{cor}

\bdem[Proof of Corollary \ref{Cor inner inner}]
Since $\vartheta(0)=0$, then $\ds\widetilde{\vartheta}=\frac{\vartheta}{z}\in \hsd$ too. By Littlewood's theorem, the composition operator 
$$
C_\fii (f)=f\circ \fii
$$ 
is a bounded operador from $\hsd$ into itself (see for instance \cite{Eero}). Hence, $\widetilde{\vartheta}\circ \fii\in \hsd$. Since $\psi = \fii (\widetilde\vartheta \circ \fii)$, we conclude that 
$\fii| \psi$. So, by the Proposition \ref{medidas en M(K)}, $\mu^\psi_1 \in M(k)$.  Taking $\alpha \psi=(\alpha \vartheta)\circ\fii$ instead of $\psi$, we get the result for the other Clark measures.
\edem

\medskip

\begin{rem}
We can get the conclusion of the previous proposition using the Clark measures of the composition. This point of view is related with the last results in \cite{H-J-W}.  
Let $\fii$ and $\vartheta$ be inner functions such that $\vartheta(0)=0$. As we have already mentioned, the composition $\psi = \vartheta \circ \fii $ is an inner function. On the other hand the Clark measures 
$\lel \mu_\al^\psi \ril_{\al \in \T}$ are given by the following formula 
$$
\mu_\al^\psi = \int\limits_\T d\mu_\beta^\fii \, d\mu_\al^\vartheta (\beta).
$$ 
Indeed, for the different values of  $\al \in \T$, let $\mu_\al^\vartheta$, $\mu_\al^\psi$, and $\mu_\al^\fii$ denote the Clark measures of $\vartheta$, $\psi$, and $\fii$, respectively. Then
\begin{align*}
\int_\T \frac{1-\modu{z}^2}{\modu{\z -z}^2} \, d\mu_\al^\psi(\z)
&= \frac{1- \modu{\psi(z)}^2}{\modu{\al - \psi(z)}^2}
= \frac{1 - \modu{\vartheta \circ \fii (z)}^2}{\modu{\al - \vartheta \circ \fii (z)}^2}\\
&= \frac{1 - \modu{\vartheta(\fii(z))}^2}{\modu{\al - \vartheta(\fii (z))}^2}
= \int_\T \frac{1-\modu{\fii(z)}^2}{\modu{\beta - \fii(z)}^2} \, d\mu_\al^\vartheta(\beta)\\
&= \int_\T \int_\T \frac{1-\modu{z}^2}{\modu{\z - z}^2} \, d\mu_\beta^\fii(\z) \, d\mu_\al^\vartheta(\beta).
\end{align*}  
Hence,
$$
\int_\T P_z(\z) \, d\mu_\al^\psi(\z) = \int_\T \int_\T P_z(\z) \, d\mu_\beta^\fii(\z)\, d\mu_\al^\vartheta(\beta).
$$ 
We obtain the result because the family $\lel P_z \ril_{z\in \D}$ dense in $C(\T)$.
Since $\vartheta(0)=0$, its Clark measures are probability measure. Therefore, we get that  
$$
\mu_\al^\psi \in \overline{\text{c.c} \left(\lel \mu_\beta^\fii\ril_{\beta\in \T}\right)}.
$$ 
By Theorem \ref{eq. con Iso}, each Clark measure $\mu_\beta^\fii$ is included in $M(k^\fii)=\text{Iso}(\fii)$. Since $\text{Iso}(\fii)$ is convex and $w^*$-closed (see \cite{Alek} or Remark \ref{obs. Teo Alek} below), it holds that $\mu_\al^\psi \in M(k^\fii)\subset M(k)$ for every $\al \in \T$.   
\end{rem}

\subsubsection{Proof of Theorem  \ref{eq. con Iso}}

Recall that given a Borel (complex) measure on $\T$ we define the Cauchy transform of $\mu$ as
$$
K(\mu)(z)=\int_\T \frac{1}{1-z\overline{\z}}\,d\z,
$$
where $z\in\D$. By a well know theorem due to Smirnov, $K(\mu)\in H^p$ for every $p\in (0,1)$ (see for instance \cite[Thm 4.2.2]{Nik1}, 
\cite[Thm 2.1.10]{C-M-R}). We will say that the inner function $\fii$ is a divisor of a function $f\in H^p$ if $f \fii^{-1}\in H^p$ for some 
$p$. Note that this is compatible with the notion of division between inner functions used before. The following result was proved by Aleksandrov in \cite{Alek}.

\begin{pro}\label{lemas 1 y 2}
Let $\mu$ be a (complex) Borel measure on $\T$ and let $\fii$ be an inner function. If $\fii$ is a divisor of $K(\mu)$ then every function in $\model{\fii}\cap H^\infty$ has non-tangential limits $|\mu|$-almost everywhere on $\T$ and for this angular boundary value we have that
$$
\int_\T \overline{f(\z)}\, d\mu(\z)=0.
$$
\end{pro}

\medskip

Recall that $H^\infty$ consists on the Banach algebra of bounded holomorphic functions defined in $\D$. Note that, in particular, as a consequence of the previous Proposition, all the reproducing kernels have non-tangential limits $|\mu|$-almost everywhere on $\T$. Using this result we can prove the following technical lemma.

\begin{lem}\label{lemita}
Let $\mu\in M(k^\fii)$. Then, there exists $\fii^*\in L^\infty(\T,\mu)$ such that 
$$
\lim_{r\to 1^-}\|\fii^*-\fii_r\|_{L^2(\T,\mu)}=0
$$ 
as well as
$$
\lim_{\angle z \rightarrow \w} \fii(z)=\fii^*(\w),
$$
for $\mu$-almost every $\w\in\T$.
\end{lem}
 
\bdem
Take $z\in \D$ such that $\fii(z)\neq 0$. Since $\ds k^\fii_z(w)=\frac{1-\overline{\fii(z)}\fii(w)}{1-\overline{z}w}$ and
$$
\lim_{r\to 1^-}\|(k^\fii_z)^*-(k^\fii_z)_r\|_{L^2(\T,\mu)}=0,
$$
straightforward computations led to 
$$
\lim_{r\to 1^-}\|\fii^*-\fii_r\|_{L^2(\T,\mu)}=0,
$$ 
for some measurable function $\fii^*:\T\to\C$. Since there exists a sequence $\{r_n\}$ such that $\{\fii_{r_n}\}$ converges $\mu$-almost everywhere to $\fii^*$ we can also conclude that $|\fii^*(\w)|\leq 1$ for almost every $\w\in\T$. Moreover, for every $z\in\D$
$$
(k^\fii_z)^*(\w)=\frac{1-\overline{\fii(z)}\fii^*(\w)}{1-\overline{z}\w}
$$
where the identity is understood $\mu$-almost everywhere. Since $\mu\in M(k^\fii)$ we know that
$$
\int_\T \Big|\frac{1 - \overline{\fii(z)} \fii^*(\z)}{1-\overline{z} \z}\Big|^2 d\mu(\z) 
= \frac{1 - \modu{\fii(z)}^2}{1 - \modu{z}^2}.
$$
Following Aleksandrov, we define the function $\Phi:\D^2\to\C$ by 
$$
\Phi(z,w)=\int_\T \frac{1 - \overline{\fii(\overline{z})} \fii^*(\z)}{1-z \z}\,  \frac{1 - \fii(w)\overline{\fii^*(\z)}}{1-w\overline{\z} }d\mu(\z)- 
\frac{1 - \overline{\fii(\overline{z})} \fii(w)}{1-z\, w} .
$$
This function is holomorphic in the bidisk and $\Phi(\overline{z},z)=0$. Therefore, $\Phi(z,w)$ is identically equal to zero in $\D^2$, in particular 
\begin{align*}
0&=\overline{\Phi(\overline{z},0)}=\int_\T \frac{1 - \fii(z) \overline{\fii^*(\z)}}{1-z \overline{\z}}\,  \big(1 - \overline{\fii(0)}\fii^*(\z)\big)d\mu(\z)- 
\big(1 - \fii(z) \overline{\fii(0)}\big).
\end{align*}
If $\mu=\mu_s+\mu_{ac}$ is the Lebesgue decomposition of $\mu$ with respect to the Lebesgue measure $m$ on $\T$, and we define 
$$
\nu=\big(1 - \overline{\fii(0)}\fii^*\big)\mu-\chi_{\text{supp}(\mu_s)^c}\,m,
$$ 
then the previous identity can be written as
\begin{align}
\int_\T \frac{1 - \fii(z) \overline{\fii^*(\z)}}{1-z \overline{\z}}\,d\nu(\z)=0.\label{eq para Cauchy}
\end{align}
Note that the singular parts of $\mu$ and $\nu$ are mutually absolutely continuous. As we mention in the section of preliminaries, we already know that 
$$
\lim_{\angle z \rightarrow \w} \fii(z)=\fii^*(\w)
$$
holds $m$-almost every $\w\in\T$. So, it only remains to prove that the non-tangential limit holds $|\nu|$-almost everywhere, because this implies that the limit holds $\mu_s$-almost everywhere too. With this aim, note that the identity\eqref{eq para Cauchy} can be rewritten as
$$
K(\nu)(z)=\fii(z) K(\fii^*\nu).
$$
Then, we can apply Proposition \ref{lemas 1 y 2} to the measure $\nu$, and we get that for every $z\in\D$ 
$$
\lim_{\angle w \rightarrow \z} k^\fii_z(w)=(k^\fii_z)^*(\z),
$$
and this directly implies that the limit
$
\ds \lim_{\angle z \rightarrow \w} \fii(z)=\fii^*(\w)
$
also holds.\edem

\medskip
As a consequence of this lemma, more precisely as a consequence of its proof, we have the following corollary.

\begin{cor}\label{cor del lemita}
Let $\mu\in M(k^\fii)$. Given any $z\in \D$
$$
(k^\fii_z)^*(\w)=\frac{1-\overline{\fii(z)}\fii^*(\w)}{1-\overline{z}\w},
$$
for $\mu$-almost every $\w\in\T$.
\end{cor} 

To conclude the proof of Theorem \ref{eq. con Iso}, we also need the following result due to Aleksandrov (see \cite{Alek})

\begin{teo}[Aleksandrov  \cite{Alek}]\label{Teo. de Aleksandrov}
If $\mu$ is a positive measure on $\T$ and $\fii$ is an inner function then the following three assertions are equivalent :
\begin{enumerate}
\item[i.)] $\mu \in \text{Iso}(\fii)$.
\item[ii.)] The function $\fii$ has angular boundary values $\mu$-almost everywhere on $\T$ and for all $z\in \D$
$$
\int_\T \modu{\frac{1 - \overline{\fii(z)} \fii(\z)}{1-\overline{z} \z}}^2 d\mu(\z) 
= \frac{1 - \modu{\fii(z)}^2}{1 - \modu{z}^2}.
$$
\item[iii.)] There exists a function $\theta \in H^\infty$ such that $\norm{\theta}_{H^\infty} \leq 1$ and for all $z \in \D$ 
$$
\int_\T P_z(\z) \, d\mu(\z) = \Preal \left( \frac{1 + \theta(z) \fii(z)}{1-\theta(z) \fii(z)} \right).
$$
\end{enumerate}
\end{teo}

\medskip

\begin{rem}\label{obs. Teo Alek}
Note that by the item (iii), if $\mu \in \text{Iso}(\fii)$ then 
$$
\norm{\mu} = \Preal \left( \frac{1 + \theta(0) \fii(0)}{1-\theta(0) \fii(0)} \right) 
=\frac{1- \modu{\theta(0) \fii(0)}^2}{\modu{1 -\theta(0) \fii(0)}^2}.
$$
So, $\text{Iso}(\fii)$ is included in the ball with radius $R=(1 -|\fii(0)|)^{-2}$. In particular, it is $w^*$-compact.
\end{rem}

\medskip

In Theorem \ref{eq. con Iso}, the proof that $\norm{\Phi(f)-f_r}_{L^2(\T,\mu)} \xrightarrow[r\rightarrow 1^-]{} 0 $ for every $f\in \model{\fii}$, when $\mu\in  \text{Iso}(\fii)$, follows essentially the lines described in \cite{Poltoratskii} for Clark measures. However, for the convenience of the reader, we include a detailed proof adapted to our setting. Firstly, we prove the next technical result.

\begin{lem}\label{lemita tecnico a lo Polto}
Let $g\in \model{\fii}$, and let $g_n$ the  $n$-th partial sum of its Taylor expansion in $\D$. Then, given $\mu\in \text{Iso}(\fii)$
$$
 \lim_{n\to\infty}\|\Phi(g)-g_n\|_{L^2(\T,\mu)}=0.
$$
\end{lem}
\bdem
Recall that, if $S^*$ denotes the backward shift operator, then $S^* g\in \model{\fii}$. Thus, by Theorem \ref{Teo Alek}, for $\mu$-almost every $\w\in\T$ 
\begin{equation}\label{eq Polto}
g_n(\w) = \Phi(g)(\w) - \w^{n+1} \Phi((S^*)^{n+1}g)(\w).
\end{equation}
Therefore 
\begin{align*}
\norm{\Phi(g) - g_n}_{L^2(\mu)} = \norm{\Phi ({S^*}^{n+1} g)}_{L^2(\mu)} 
= \norm{{S^*}^{n+1} g}_{H^2} \conv 0. 
\end{align*} 
\edem

\bdem[Proof of Theorem \ref{eq. con Iso}] Given $\mu\in \text{Iso}(\fii)$, let  
$g$ be an holomorphic function defined on $\D$ and suppose that 
$$
g(z) = \sum_{k=0}^\infty c_k z^k
$$
is the  Taylor expansion of $g$ in $\D$. Applying the Abel summation method (summation by parts) to the power series of $g_r$ we get 
\begin{align*}
\norm{g_r}^2_{L^2(\mu)} 
&= \int_\T \modu{\sum_{k=0}^\infty \left(\sum_{n=0}^k c_n \z^n \right) 
\left( r^k - r^{k+1} \right) }^2 d\mu(\z).
\end{align*}
Note that $\ds{\sum_{k=0}^\infty}\left(r^k - r^{k+1} \right) =1$. So, using Jensen's inequality we obtain
\begin{align*}
\norm{g_r}^2_{L^2(\mu)}  
&\leq\int_\T \sum_{k=0}^\infty \modu{\sum_{n=0}^k c_n \z^n}^2
\left( r^k - r^{k+1} \right)  d\mu(\z) \\
&= \sum_{k=0}^\infty\int_\T  \modu{\sum_{n=0}^k c_n \z^n }^2 d\mu(\z) \left(r^k - r^{k+1} \right)\\ 
&= \sum_{k=0}^\infty \norm{g_k}^2_{L^2(\mu)} \left(r^k - r^{k+1} \right).
\end{align*}

Let $\eps>0$. By Lemma \ref{lemita tecnico a lo Polto}, there exists $n\geq 1$ such that
\begin{align*}
\norm{\Phi(f) - f_n}_{L^2(\mu)} < \frac{\eps}{4},
\end{align*} 
where  $f_n$ denotes the $n$-th partial sum of the Taylor expansion in $\D$ of the function $f$. Moreover, we can take this $n$ so that
$$
\|\Phi((S^*)^{n+1}f)\|_{L^2(\mu)}<\frac{\eps}{2}
$$
also holds. On the other hand, since the measure $\mu$ is finite,  there exists $r_0 \in (0,1)$ such that 
$$
\norm{f_n - (f_n)_r}_{L^2(\mu)} < \frac{\eps}{4},
$$ 
for every $r\in (r_0,1)$. So, we have that
\begin{align*}
\norm{\Phi(f) - f_r}_{L^2(\mu)}&\leq \frac{\eps}{2}+\|f_r-(f_n)_r\|_{L^2(\mu)}
=\frac{\eps}{2}+\|(f-f_n)_r\|_{L^2(\mu)}\\
&=\frac{\eps}{2}+ \left(\sum_{k=0}^\infty \norm{(f-f_n)_k}^2_{L^2(\mu)} \left(r^k - r^{k+1} \right)\right)^{1/2}.
\end{align*}

Using \eqref{eq Polto}, straightforward computation show that 
$$
\lim_{k\to\infty} \norm{(f-f_n)_k}^2_{L^2(\mu)}=\|\Phi((S^*)^{n+1}f)\|^2_{L^2(\mu)}.
$$ 
Therefore
$$
\limsup_{r \rightarrow 1^-} \sum_{k=0}^\infty \norm{(f-f_n)_k}^2_{L^2(\mu)} \left(r^k - r^{k+1} \right)
=\|\Phi((S^*)^{n+1}f)\|^2_{L^2(\mu)}.
$$ 
Consequently
\begin{align*}
\limsup_{r \rightarrow 1^-}\norm{\Phi(f) - f_r}_{L^2(\mu)}
&\leq \frac{\eps}{2}+ \|\Phi((S^*)^{n+1}f)\|_{L^2(\mu)}\leq \eps.
\end{align*}


In particular this holds for the reproducing kernels. On the other hand, since the spaces  $\model{\fii}$ and $L^2(\T,\mu)$ are isometrically isomorphic, the kernels clearly reproduce themseves with respect to $\mu$.

\medskip

Conversely, take $\mu\in M(k^\fii)$. By Lemma \ref{lemita} and Corollary \ref{cor del lemita}, we get that the condition (ii) of Theorem \ref{Teo. de Aleksandrov} holds. Then $\mu\in \text{Iso}(\fii)$.
\edem


\begin{center}
References
\end{center}

\fontsize{10}{10}\selectfont \small{
}

\end{document}